\documentclass{article}

\usepackage{amsmath, amsthm, amsfonts}
\usepackage{color}

\usepackage{hyperref}

\usepackage{longtable}

\usepackage{csquotes}

\usepackage{mdframed}
\theoremstyle{definition}
\newmdtheoremenv{boxProb}{Problem}
\newmdtheoremenv{boxDef}{Definition}
\newmdtheoremenv{boxCor}{Corollary}
\newmdtheoremenv{boxThm}{Theorem}
\newmdtheoremenv{compjob}{Computational Job}
\newmdtheoremenv{reqi}{Requirement}

\usepackage{graphicx}
\usepackage{subfigure}          

\usepackage{placeins}

\usepackage{tabularx}

\usepackage{lmodern,textcomp}

\usepackage{algorithm}
\usepackage{algpseudocode}

\usepackage{enumitem}

\usepackage{xspace} 
\newcommand\largeparbreak{\par\bigskip}

\newcommand{\inv}{^{-1}\xspace}

\newcommand{\bmuL}{{\boldsymbol{\mu}_{\textsf{L}}}\xspace}
\newcommand{\bmuR}{{\boldsymbol{\mu}_{\textsf{R}}}\xspace}

\newcommand{\blambda}{\boldsymbol{\lambda}\xspace}
\newcommand{\bxi}{\boldsymbol{\xi}\xspace}

\renewcommand{\t}{^\textsf{T}\xspace}

\newcommand{\hblambda}{\hat{\blambda}\xspace}

\newcommand{\away}[1]{}

\newcommand{\R}{\mathbb{R}\xspace}

\newcommand{\cL}{\mathcal{L}\xspace}

\newcommand{\cF}{\mathcal{F}\xspace}

\newcommand{\cO}{\mathcal{O}\xspace}

\newcommand{\bS}{\textbf{S}\xspace}

\newcommand{\bx}{\textbf{x}\xspace}
\newcommand{\bxL}{\textbf{x}_{\textsf{L}}\xspace}
\newcommand{\bxR}{\textbf{x}_{\textsf{R}}\xspace}

\newcommand{\br}{\textbf{r}\xspace}

\newcommand{\bp}{\textbf{p}\xspace}

\newcommand{\bz}{\textbf{z}\xspace}
\newcommand{\be}{\textbf{1}\xspace}
\newcommand{\bei}[1]{{\textbf{e}}\xspace}
\newcommand{\bw}{\textbf{w}\xspace}

\newcommand{\bH}{\textbf{H}\xspace}
\newcommand{\bI}{\textbf{I}\xspace}

\newcommand{\bO}{\textbf{0}\xspace}

\newcommand{\opdiag}{\textsl{diag}\xspace}

\newcommand{\tol}{{\textsf{tol}}\xspace}

\newcommand{\ttau}{\tilde{\tau}\xspace}
\newcommand{\tomega}{\tilde{\omega}\xspace}
\newcommand{\trho}{\tilde{\rho}\xspace}

\newcommand{\PBNLP}{\textsf{PBNLP}}
\newcommand{\NLP}{\textsf{NLP}}

\usepackage{todonotes}

\title{Interior Point Method with Modified Augmented Lagrangian for Penalty-Barrier Nonlinear Programming}

\author{Martin Neuenhofen}

\begin{document}

\maketitle

\begin{abstract}
We present a numerical method for the local solution of nonlinear programming problems. The SUMT approach of Fiacco and McCormick results in a merit function with quadratic penalties and logarithmic barriers. Our NLP solver works by directly minimizing this merit function.

In our method, we use different concepts that each shall aim at the efficient treatment of one respective special feature of this merit function. The features are: large quadratic penalty terms, and badly scaled logarithmic barriers.

The quadratic penalties are treated with a modified Augmented Lagrangian technique. It enables large step sizes despite nonlinearity of the equality constraints. The logarithmic barriers we treat with a primal-dual interior-point path-following technique.

We prove global convergence of the method and local quadratic convergence. We further prove weak polynomial time-complexity in the special case where the nonlinear program is a linear program.

We also use a trust-funnel so to avoid that the method converges to any stationary points that are infeasible to the constraints.
\end{abstract}


\section{Introduction}

\paragraph{Problem statement}
This paper presents a method for the solution of the penalty-barrier nonlinear programming problem
\begin{equation}
\tag{\PBNLP}\label{eqn:PBNLP}
\begin{aligned}
	&\operatornamewithlimits{min}_{\bx \in \overline{\Omega}}& \phi(\bx) := &\quad f(\bx) + \frac{\rho}{2} \cdot \|\bx\|_\bS^2 \\
	& & &\quad + \frac{1}{2 \cdot \omega} \cdot \|c(\bx)\|_2^2 \\
	& & &\quad - \tau_E \cdot \be\t \cdot \Big(\,\log(\bx-\bxL) - \log(\bxR - \bx)\,\Big)\,,
\end{aligned}
\end{equation}
where
$$ \Omega := \lbrace\,\bxi \in \R^n \, \vert \, \bxL < \bxi < \bxR\, \rbrace\,. $$
In this problem statement, we are provided with: twice smooth differentiable functions $f\,:\,\R^n \rightarrow \R$, $c\,:\,\R^n \rightarrow \R^m$, and $\cL(\bx,\blambda) := f(\bx)-\blambda\t\cdot c(\bx)$\,; a matrix $\bS \in \R^{n \times n}$ symmetric positive definite (spd), and small parameters $\rho,\omega,\tau_E>0$\,, e.g. $\rho=\omega=\tau_E=10^{-8}$\,; an initial guess $\bx_0 \in \Omega$ such that $$\|c(\bx_0)\|_\infty \leq 2 \cdot \epsilon$$
for some user-specified funnel-parameter $\epsilon>0$\,, e.g. $\epsilon=0.1$\,.

\paragraph{Motivation}
As motivated in \cite{PPDIPM}, the minimization of \eqref{eqn:PBNLP} is related to locally minimizing the nonlinear program
\begin{equation}
\tag{\NLP}\label{eqn:NLP}
	\begin{aligned}
		&\operatornamewithlimits{min}_{\bx \in \R^n}& f(\bx)& \\
		&\text{s.t.} 								& c(\bx)&=\bO\,,\\
		& 											& \bxL \leq \phantom{c(}\bx\phantom{)}&\leq \bxR\,.
	\end{aligned}
\end{equation}
But in contrast to \eqref{eqn:NLP}, the minimization of \eqref{eqn:PBNLP} is also meaningful when $m>n$ or when $c$ has no root. Hence, minimizing \eqref{eqn:PBNLP} is more generic in approach than solving \eqref{eqn:NLP}. Also, it circumvents technicalities related to constraint qualifications, local uniqueness of a local minimizer for $\bx$, and issues related to degeneracy.

\paragraph{Challenges}
As on the contra side, according to the author, the minimization of \eqref{eqn:PBNLP} can be more challenging than the solution of the \eqref{eqn:NLP}. As discussed in \cite{PPDIPM}, this is because there is no such thing as a projector that could move a trial point back into the domain of feasible iterates with respect to $c$. Hence, nonlinearities in $c$ can make it challenging to achieve progress in the line-search.
As another issue, $\phi$ can have local minimizers that do not yield small values for $\|c(\bx)\|_2$\,. Approaches must be utilized that hinder the iterates from converging into such infeasible stationary points.
Further, the logarithmic barriers demand for a special treatment. This is because they lead to unbounded values when the iterates approach $\partial\Omega$.
As a final issue, the small parameters $\omega,\tau_E$ result in bad scaling of $\phi$.

\paragraph{Background}
The formulation of a function like $\phi$ and the relation between \eqref{eqn:PBNLP} and \eqref{eqn:NLP} is analyzed extensively in \cite{SUMT}. However, their analysis is primarily for the case where $m<n$, which appears unnecessarily restrictive when only considering \eqref{eqn:PBNLP}.

In \cite{myOCP} we make use of formulations of \eqref{eqn:PBNLP} with $m>n$ to solve optimal control problems. Larger values of $m$ result from higher orders of quadrature for continuous equality constraints. These higher orders are desirable to improve accuracy.

\largeparbreak

In the quest of finding or deriving efficient solution methods for \eqref{eqn:PBNLP}, it seems reasonable to first attempt using those methods from the literature that solve \eqref{eqn:NLP} by actually minimizing \eqref{eqn:PBNLP}. Examples of such methods are given in \cite{ForsgrenGill,ChenGoldfarb,myCQP}. All these methods are interior-point methods. In attempt to also be able to apply techniques from successive quadratic programming \cite{SNOPT}, we developed the method in \cite{PPDIPM}. This method minimizes $\phi$ in a direct way, using a line-search technique and step-directions computed as solutions of a convex quadratic program.

\largeparbreak

Two technical details keep it difficult to treat \eqref{eqn:PBNLP}, of which one is inherited from \eqref{eqn:NLP}.

First, since $m>n$ and $\omega$ is small, it is difficult to treat nonlinearities in $c$ within the line-search. In particular, when $\omega$ is small then the line-search will only accept small step sizes. Conventionally, this would be overcome by using a second-order correction technique \cite{NumOpt,TrustRegionMethods}. But since $m>n$, there is no way to project back an infeasible iterate onto a feasible sub-domain. Even further, since $c$ is treated with an inexact penalty, even a solution of $c(\bx)=\bO$ existed, it is unlikely that the minimizer of $\phi$ yields $c(\bx)=\bO$. The Augmented Lagrangian technique is an elegant alternative to quadratic penalties, at it often allows for keeping the penalty terms small \cite{ALM_IPM1,Lancelot}. But, the Augmented Lagrangian technique is an exact penalization technique, whereas in \eqref{eqn:PBNLP} we deliberately use an inexact penalty. In \cite{MALM} we generalized the Augmented Lagrangian technique to this purpose. We use it in the here proposed method because we hope that it allows for better progress during the line-search.

The second issue, that is inherited from \eqref{eqn:NLP}, is with second-order local convergence. Interior-point methods base on second-order local convergence of the Newton iteration when solving the barrier-KKT conditions for a respective value of $\tau$, where $\tau \searrow \tau_E$. This has drawbacks, as we rely on efficient ways to decrease $\tau$. And, suitable initial values for $\tau$ in conjunction with suitable sufficiently central initial guesses for $\bx_0$ must be found (which yet still is an open question for \eqref{eqn:NLP}). In contrast, when using a successive quadratic programming method (SQP)\footnote{SQP, not Actice Set Method; i.e. in our understanding of SQP there are no active sets involved! Active set methods are unacceptable for the large-scale problems that we are interested in treating. This is because in the average case it takes to many iterations for the active set to converge, namely in practice one often observes it to take $\cO(n)$ iterations \cite{NumOpt}.} we are not aware of any existing method that has polynomial time-complexity per iteration and achieves a locally fast rate of convergence; where with fast we mean a rate of convergence that is at least linear with a rate of contraction of $0.5$\,, i.e. invariant to conditioning or scaling of the problem instance. We asked a lot of leading experts on this, we offered considerable money for a solution, and we wrote a research proficiency examination report of 27 pages with 74 references on this matter, which yet remains to publish. In result, we honestly believe that such an SQP method cannot exist. Hence, here we focus on an interior-point method, despite all the drawbacks of using interior-points.

\paragraph{Outline}
In this paper we present a method that minimizes \eqref{eqn:PBNLP}. Our method uses two techniques: the augmented Lagrangian technique \cite{Renke,ALM_IPM1} and the primal-dual path-following technique \cite{AdlerMonteiro}. Globalization is achieved through a line-search with a merit function, cf. \cite{TrustRegionMethods}.

Without further emphasize, we also included a feasibility-funnel that hinders the iterates from moving away too far from approximate minimizers of $\|c(\bx)\|_2$. This is achieved with a mildly-scaled (parameter $\ttau > \tau$) primal barrier term. In the convergence theory, it could be interpreted as a part of $f$, which is why we just include it in our formulas for sake of easy implementability, but ignore it completely within our further discussion and convergence analysis. Since mildly scaled, it should not spoil the convergence. And since minimization of $\phi$ aims at small values for $c(\bx)$ anyway, the additional barrier terms should not have a big influence on the actual minimizer $\bx$ of $\phi$.

\section{The method}

\paragraph{Preliminaries}
Usually, interior-point methods are nested iterative schemes of two nested levels \cite{ChenGoldfarb}:
\begin{itemize}
	\item Outermost level: On the outermost iteration level an update the barrier parameter $\tau$ is performed.
	\item Inner level: Within that outer iteration there is an inner iteration that iteratively solves the optimization problem for the current barrier parameter.
\end{itemize}
In contrast to that, in this paper we have a nested iterative scheme of three levels:
\begin{itemize}
	\item Outermost level: On the outermost iteration we update the barrier parameter.
	\item Outer level: Within the outermost iteration there is the outer iteration. It performs updates of the augmented Lagrangian multipliers.
	\item Inner iteration: Within the outer iteration there is the inner iteration. It solves the optimization problem for the current values of barrier parameter and augmented Lagrangian multipliers.
\end{itemize}

We consider three types of numbers:
\begin{itemize}
	\item Constants: These numbers never change their values. $\epsilon>0$ is the width of the feasibility funnel. $\tau_E>0$ is the final barrier parameter. $\tol>0$ is the tolerance for the KKT residual. $\rho>0$ is the convexization parameter. $\omega>0$ is the quadratic penalty parameter. $\tomega > \omega$ is the relaxed penalty parameter for the Augmented Lagrangian framework. $\ttau$ is the relaxed barrier parameter for the feasibility funnel. $\nu>0$ is the dual weighting of the primal-dual merit function that our method uses. $\vartheta \in (0,1)$ is the fraction-to-the-boundary coefficient for the line-search.
	\item Parameters: These numbers only change during the outer or outermost iterations. $\tau\geq\tau_E$ is the current barrier-parameter. It is geometrically decreased during the outermost iteration, while inner iterations attempt to find solutions for $\bx$ on the central path for this value of $\tau$. $\hblambda \in \R^m$ is the current vector of augmented Lagrangian multipliers. It is updated during outer iterations. Again, the inner iterations will try to find solutions for $\bx$ on the central path for this value of $\hblambda$.
	\item Variables: These numbers are updated during every inner iteration. $\bx \in \Omega$ is the current iterate for a solution of \eqref{eqn:PBNLP}. $\blambda \in \R^m$ are the $\tomega$-penalty regularized Lagrange multipliers for the equality constraints. $\bmuL,\bmuR>\bO$ are the $\tau$-barrier regularized Lagrange multipliers for the box-constraints $\bxL \leq \bx \leq \bxR$.
\end{itemize}

We consider primal-dual iterates:
\begin{align*}
	\bz = \begin{pmatrix}
		\bx\\
		\blambda\\
		\bmuL\\
		\bmuR
	\end{pmatrix} \in \R^{n + m + n + n}
\end{align*}
More particular, during all iterations it will hold $\bx \in \Omega$ and $\bmuL,\bmuR>\bO$. The set of iterates $\bz$ obeying to these conditions is called $\cF \subset \R^{3\cdot n + m}$.

We also consider a parameter vector:
\begin{align*}
	\bp = \begin{pmatrix}
		\tau\\
		\hblambda
	\end{pmatrix} \in \R^{1+m}
\end{align*}
For a compact notation we define
\begin{align*}
	\bw(\bx) := \frac{\ttau}{\epsilon\cdot\be+c(\bx)}-\frac{\ttau}{\epsilon\cdot\be-c(\bx)}\,.
\end{align*}

We use this root function:
\begin{align*}
	F(\bz,\bp) := \begin{pmatrix}
		\nabla_\bx \cL\big(\,\bx\,,\,\hblambda+\blambda + \bw(\bx) \,\big) + \rho \cdot \bS \cdot \bx - \bmuL + \bmuR\\
		c(\bx) + \omega \cdot \hblambda + (\omega + \tomega) \cdot \blambda\\
		\bmuL \cdot (\bx-\bxL) - \tau \cdot \be\\
		\bmuR \cdot (\bxR - \bx) - \tau \cdot \be
	\end{pmatrix} =: \begin{pmatrix}
		\br_{\textsf{Dual}}\\
		\br_{\textsf{Primal}}\\
		\br_{\textsf{Comp,L}}\\
		\br_{\textsf{Comp,R}}
	\end{pmatrix}
\end{align*}
Using the writing $\bH := \nabla_\bx^2 \cL\big(\,\bx\,,\,\hblambda+\blambda+\bw(\bx)\,\big)$, it has the following Jacobian with respect to $\bz$:
\begin{align*}
	DF(\bz,\bp) = \left[\begin{array}{c|ccc}
		\bH + \rho \cdot \bS & -\nabla c(\bx) & -\bI & \bI\\
		\hline
		\nabla c(\bx)\t & -(\omega + \tomega) \cdot \bI & \bO & \bO\\
		\opdiag(\bmuL) & \bO & \opdiag(\bx-\bxL) & \bO\\
		-\opdiag(\bmuR) & \bO & \bO & \opdiag(\bxR-\bx)
	\end{array}\right]
\end{align*}

We define the following primal-dual merit function:
\begin{align*}
M(\bz,\bp) :=
	&\quad \cL(\bx,\hblambda)\\
	&\quad -\ttau \cdot \be\t \cdot \Big(\,\log\big(\epsilon \cdot \be + c(\bx)\big) + \log\big(\epsilon\cdot\be - c(\bx)\big)\,\Big)\\
	&\quad + \frac{1}{2 \cdot \tomega} \cdot \big\|\,c(\bx) + \omega \cdot (\hblambda + \blambda)\,\big\|_2^2\\
	&\quad + \frac{\rho}{2} \cdot \|\bx\|_\bS^2 + \frac{\omega}{2} \cdot \|\blambda\|_2^2 - \tau \cdot  \be\t \cdot \Big(\,\log(\bx-\bxL) + \log(\bxR-\bx)\,\Big)\\
	&\quad + \frac{\nu}{2 \cdot \tomega} \cdot \big\|\,c(\bx) + \omega \cdot \hblambda + (\omega + \tomega) \cdot \blambda\,\big\|_2^2\\
	&\quad - \nu \cdot \tau \cdot \be\t \cdot \Big(\,\log\big(\bmuL \cdot \frac{\bx-\bxL}{\tau}\big)+\be-\bmuL \cdot \frac{\bx-\bxL}{\tau}\,\Big)\\
	&\quad - \nu \cdot \tau \cdot \be\t \cdot \Big(\,\log\big(\bmuR \cdot \frac{\bxR-\bx}{\tau}\big)+\be-\bmuR \cdot \frac{\bxR-\bx}{\tau}\,\Big)
\end{align*}

\paragraph{Iterative approach}

The three iterative levels are given in the code framework below: the outermost, the outer, and the inner iteration.
\begin{algorithmic}[0]
	\Procedure{ProposedMethod}{}
	\State Given: initial guess $\bz \in \cF$, parameters $\bp$, constants $\tau_E,\tol$.
	\While{$\tau>\tau_E$}
		\State Outermost iteration: Update $\tau$.
		\While{$\|\blambda\|_\infty>\tol$}
			\State Outer iteration: Update $\hblambda$.
			\While{$\|F(\bz,\bp)\|_\infty>\tol$}
				\State Inner iteration: Update $\bz$.
			\EndWhile
		\EndWhile
	\EndWhile
	\State \Return $\bx$
	\EndProcedure
\end{algorithmic}

Each iteration is embedded in a while-loop, whose condition to achieve is the scope of the respective update. In accurate terms: The outermost iteration updates $\tau$ such that eventually $\tau = \tau_E$. The outer iteration updates the Augmented Lagrangian $\hblambda$ such that the root $\bz$ of $F(\cdot,\bp)$ will eventually have a component $\blambda$ that approaches zero. The inner iteration computes $\bz$ with a globalized Quasi-Newton iteration, such that $F(\bz,\bp)$ approaches $\bO$ for the given parameters $\bp$. When all the three while-conditions are satisfied, then the component $\bx$ of $\bz$ is an accurate minimizer of \eqref{eqn:PBNLP}.

In the following we describe the iterations in detail.

\subsection{Inner iteration}
Given a parameter vector $\bp$ and an initial guess $\bz \in \cF$, we use a globalized Quasi-Newton iteration to find a root $\bz \in \cF$ of $F(\cdot,\bp)$. The globalization works through a line-search with the merit function $M$.

In a conventional Newton iteration we would compute an update $\Delta\bz$ as the solution of
$$ DF(\bz,\bp) \cdot \Delta\bz = -F(\bz,\bp)\,.$$
However, the iteration $ \bz := \bz + \Delta\bz$ may not be globally convergent, since the Newton iteration is not globally convergent.

To achieve global convergence, we employ a concept from \cite{ForsgrenGill}. It works by reformulating the Newton system for $\Delta\bz$. Rescaling the latter two rows and changing the signs in the second and third columns, and adding a shift $(\trho-\rho) \cdot \bS$ to the upper left block, we arrive at the following system:
\begin{align*}
	&\left[\begin{array}{c|ccc}
	\bH + \trho \cdot \bS & -\nabla c(\bx) & \bI & \bI\\
	\hline
	\nabla c(\bx)\t & -(\omega + \tomega) \cdot \bI & \bO & \bO\\
	\bI & \bO & -\opdiag\Big(\frac{\bx-\bxL}{\bmuL}\Big) & \bO\\
	\bI & \bO & \bO & -\opdiag\Big(\frac{\bxR-\bx}{\bmuR}\Big)
	\end{array}\right] \cdot \begin{pmatrix}
		\Delta\bx\\
		-\Delta\blambda\\
		-\Delta\bmuL\\
		\Delta\bmuR
	\end{pmatrix}\\
	 = &-\begin{pmatrix}
		\br_{\text{Dual}}\\
		\br_{\text{Primal}}\\
		\phantom{-}\opdiag(\bmuL)\inv\cdot\br_{\text{Comp,L}}\\
		-\opdiag(\bmuR)\inv\cdot\br_{\text{Comp,R}}		
	\end{pmatrix}
\end{align*}
As we see from the saddle structure of the matrix, the system has at least $2 \cdot n + m$ strictly negative eigenvalues. We choose the variable $\trho\geq \rho$ large enough such that the remaining $n$ eigenvalues are all strictly positive. In \cite{ForsgrenGill,MALM} it is shown that then $\Delta\bz$ is a descent direction for the merit function $M$.

Our inner iteration works by computing $\Delta\bz$ from the above system for a sufficiently large value of $\trho$\,. Afterwards, we perform a back-tracking line-search along the merit function $M$\,. We subject the step-size to the Armijo condition. The initial trial step size is $\alpha_0 = \min(1,\vartheta \cdot \alpha_{\textsf{Box}})$, where $\alpha_{\textsf{Box}}>0$ is the maximum step size such that $\bz + \alpha_{\textsf{Box}} \cdot \Delta\bz \in \overline{\cF}$.

The inner iteration is converged when $\|F(\bz,\bp)\|_\infty\leq \tol$.

\subsection{Outer iteration}

We are given two vectors $\bz,\bp$, such that $F(\bz,\bp)\approx \bO$. The outer iteration is a modified augmented Lagrangian method (MALM) \cite{MALM}.

The goal of the MALM iteration is to successively refine the vector $\hblambda$, such that the root $\bz$ of $F(\cdot,\bp)$ satisfies $\blambda = \bO$, as holds in the limit for the conventional Augmented Lagrangian methods \cite{ALM_IPM1,Renke}.

According to \cite{ALM_IPM1}, the update would be
\begin{align*}
	\hblambda := & \hblambda + \alpha \cdot \blambda\,,\\
	\blambda := & (1-\alpha) \cdot \blambda\,,
\end{align*}
where $\alpha = 1$\,. We relax this update rule by allowing values of $\alpha \in (0,1]$\,. After the update, we have new vectors $\bz_{\textsf{new}},\bp_{\textsf{new}}$.

The concern with the update is that the new vectors may result in large values of $\|F(\bz_{\textsf{new}},\bp_{\textsf{new}})\|_\infty$. I.e., the update of the Augmented Lagrangian may completely destroy the centrality property of the iterate $\bz$. However, the centrality property is a crucial property of $\bz$ that must hold for the inner iterations to make converge at a reasonable rate. Destruction of centrality could lead to stall or impractically slow progress within the next sequence of inner iterations.

To address the concern, we make use of the relaxed update for $\hblambda$ by employing a back-tracking line search for $\alpha$. A trial value for $\alpha$ yields a trial point 
$$ \bp(\alpha) := \bp + \alpha \cdot \begin{pmatrix}
	0\\
	\blambda
\end{pmatrix}\,,$$
and analogously a trial vector for $\bz$, where only $\blambda$ is affected by the update.
Certainly, for a sufficiently small value of $\alpha>0$ the destruction of the centrality property is mild and the inner iterations can make good progress. But, choosing $\alpha$ very small has the disadvantage that the MALM iteration only converges slowly. This is why we do not just want to rely on choosing $\alpha$ small.

This is why, as a further enhancement, we compute $\bz_{\textsf{new}}$ from scratch for every trial step-size $\alpha$ by considering an auxiliary problem. That is, for every trial step size $\alpha$ we compute $\bp(\alpha)$ and then determine $\bz_{\textsf{trial}}$ as the solution of a set of equations that approximately represent the system
$$ F\big(\,\bz_{\textsf{trial}},\bp(\alpha)\,\big) = \bO\,. $$
This set of equations we write below:
\begin{subequations}
	\label{eqn:KKT_LP}
	\begin{align}
		\rho \cdot \bS \cdot \bx_{\textsf{trial}} + \nabla_\bx\cL(\bx,\hblambda_{\textsf{trial}}) - \nabla c(\bx_0) \cdot \blambda_{\textsf{trial}} - \bmuL_{\textsf{trial}} + \bmuR_{\textsf{trial}} &= \bO\\
		\nabla c(\bx)\t \cdot \bx_{\textsf{trial}} - \big(\,\nabla c(\bx)\t \cdot \bx - c(\bx) - \omega \cdot \hblambda_{\textsf{trial}} \,\big) + (\omega + \tomega) \cdot \blambda_{\textsf{trial}} &= \bO\\
		\opdiag(\bmuL_{\textsf{trial}}) \cdot (\bx_{\textsf{trial}}-\bxL) - \tau \cdot \be &= \bO\\
		\opdiag(\bmuR_{\textsf{trial}}) \cdot (\bxR-\bx_{\textsf{trial}}) - \tau \cdot \be &= \bO
	\end{align}
\end{subequations}
This system depends on the current iterate $\bz$ and on the trial value $\hblambda_{\textsf{trial}} := \hblambda + \alpha \cdot \blambda$\,. Based on these data, the above system defines a unique solution\footnote{since forming the regularized optimality system of a strictly convex quadratic program, that by construction has one unique local minimizer.} for $\bz_{\textsf{trial}}$. Easy to see, for $\alpha=0$ it holds $\bz_{\textsf{trial}}=\bz$ because then from $F(\bz,\bp)=\bO$ it follows that $\bz$ solves the above equations. The system can be solved in weak polynomial time complexity, for instance by utilizing our interior-point algorithm described in \cite{myCQP,PPDIPM}.

The trial step size $\alpha$ is accepted when
$$ \|F\big(\bz_{\textsf{trial}},\bp(\alpha)\big)\|_\infty \leq \chi \tol \,.$$

\subsection{Outermost iteration}
In the outermost iteration we are given points $\bz,\bp$, such that $\blambda\approx\bO$ and $F(\bz,\bp)\approx \bO$. The purpose of the iteration is to geometrically reduce the parameter $\tau$, such that eventually it approaches its target value $\tau \rightarrow \tau_E$\,. We know from \cite{Wright,IPM25ylater} that for problems such as linear programming and convex quadratic programming there exist geometric updates of the kind
\begin{align*}
	\tau := \max(\,\sigma \cdot \tau\,,\,\tau_E\,)\,,\quad \quad \text{where }\sigma < 1\,,
\end{align*}
such that centrality of the iterate $\bz$ can be reestablished within the next inner iteration. We opt for the same heuristic that is used in \cite{IPOPT}, i.e. we make an update with $\sigma = 0.1$\,. The update is performed every time after the MALM iteration converged to points $\bz,\bp$ that yield $\blambda\approx\bO$ and $F(\bz,\bp)\approx \bO$.

\section{Proof of Convergence: Global, locally quadratic, and complexity for linear programming}
Our claims, proved in this section, are in summary:
\begin{enumerate}
	\item If $\tomega$ is sufficiently large then our proposed method has global convergence towards a stationary point of $\phi$.
	\item Our proposed method converges locally quadratic.
	\item In the special case where $f$ and $c$ are linear functions, there are suitable values for $\sigma$ such that our proposed method has polynomial time complexity.
\end{enumerate}

\subsection{Global convergence}
We explain the global convergence of each iteration. In this course, we first discuss inner iterations, then outer iterations, and finally outermost iterations.

\largeparbreak

In \cite{ForsgrenGill} it is proven that the globalized Quasi-Newton iteration of the inner iteration converges to a root of $F$. Hence, the inner iteration converges.

The outer iteration is MALM. In \cite{MALM} we provide an analysis that shows that MALM is the conventional Augmented Lagrangian Method applied to an auxiliary system. It follows that global convergence of MALM is achieved whenever the conventional Augmented Lagrangian method converges for the auxiliary optimization problem. As discussed in \cite{MALM,ALM_IPM1,Renke} and the references therein, this holds when $\tomega$, the Augmented-Lagrangian-penalty parameter, is sufficiently large, so that the feasibility residual of the auxiliary optimization problem remains bounded. Hence, under the assumption on $\tomega$, the outer iteration converges. As thoroughly discussed in \cite{MALM}, it is at the heart of the Augmented Lagrangian method that the Augmented Lagrangian $\hblambda$ converges to the exact Lagrangian multipliers of the auxiliary optimization problem. Hence, the penalty parameter does not need to be increased boundlessly. Further, the equality constraints in the auxiliary problem (which are $c(\bx) + \omega \cdot \hblambda = \bO$) converge. Hence, $\blambda=\bO$ results in the limit.

Finally, there is the outermost iteration. Every time the outer iteration converges, as we showed above, the outermost iteration will decrease $\tau$ by a significant amount. Hence, eventually it will be $\tau=\tau_E$. Therefore, the outermost iteration converges as well.

When all the three iterations have converged, we have points $\bz,\bp$, such that $\tau=\tau_E$, $\blambda=\bO$, and
$$ F(\bz,\bp)=\bO\,. $$
This is equivalent to the condition $\nabla \phi(\bx) = \bO$\,. Hence, we just showed that our method converges globally to a stationary point of $\phi$.

As a remark, the feasibility-funnel with width $\epsilon$ can be chosen suitably small such that each stationary point of $\phi$ achieves small norms in $c$.

\subsection{Locally quadratic convergence}
The Quasi-Newton iteration in the inner iteration converges at a locally quadratic rate to a root of $F(\cdot,\bp)$\,, cf. \cite{ForsgrenGill}.

\subsection{Polynomial time-complexity for linear programming problems}
The scope of this subsection is to show that our proposed method has polynomial time-complexity when $f,c$ are linear and when choosing $\sigma = 1 - 0.1/\sqrt{2 \cdot n}$\,.

We underline that for linear programming problems there exist far more efficient iterative solution algorithms, that are both cheaper per iteration and converge in less iterations, cf. \cite{Wright,IPM25ylater,FriedlanderOrban}. The purpose of this subsection is rather that we want to give reason to believe that our proposed method inherits some sort of generic efficiency that makes it promising for efficiently solving \eqref{eqn:PBNLP}.

\largeparbreak

In the following we go through the algorithm and explain what happens. We will see that the overall time-complexity is polynomial:

In the beginning, we have $\bz \in \cF$ and a parameter vector $\bp$. Starting with the first outermost iteration, $\tau$ is updated to $\tau := \sigma \cdot \tau$. We then proceed with the first outer iteration, that updates $\hblambda$. As discussed in Section~2.2, the outer iteration also involves an update of $\bz$, namely we solve \eqref{eqn:KKT_LP} for $\bz_{\textsf{trial}}$ and update $\bz := \bz_{\textsf{trial}}$. Since $f,c$ are linear, the property 
$$ F(\bz_{\textsf{trial}},\bp)=\bO $$
follows immediately from \eqref{eqn:KKT_LP}. Hence, the step size $\alpha=1$ in the outer iteration will always be accepted.

Then we come to the inner iteration. The inner iteration immediately terminates because $F(\bz,\bp)=\bO$ already holds.

We come back to the outer iteration, that updates $\hblambda$. Since $\alpha=1$ will always be accepted and since the inner iteration immediately terminates, the outer iteration will converge rapidly.

Finally, we come back to the outermost iteration. Due to our mild choice of $\sigma$, we can use path-following results from \cite{myCQP,IPM25ylater,Wright}. These show that basically one Newton iteration is sufficient to refine a root $\bz$ of $F(\cdot,\bp)$ for the current value of $\tau$ in order to find an accurate root of $F(\cdot,\bp)$ for the updated value of $\tau$. Using the proposed value of $\sigma$, the number of outermost iterations is bounded by a weakly polynomial number \cite{myCQP}.

\largeparbreak

Since the number of all loop iterations is bounded by a polynomial and since each iteration can be computed in polynomial complexity, the entire algorithm has polynomial complexity.

\section{Conclusion}
In \cite{MALM} we merged the Augmented Lagrangian framework with the quadratic penalty method for minimizing functions with large quadratic penalty terms. The resulting approach we called Modified Augmented Lagrangien method (MALM).

In \cite{myCQP} we presented a stabilized and regularized primal-dual path-following framework, that --when applied to convex quadratic programs -- we showed to be polynomially efficient and numerically stable.

In \cite{ForsgrenGill} Forsgren and Gill provide a merit function that, together with an inertia correction, decreases monotonously in Quasi-Newton directions computed from our stabilized and regularized path-following equations.

In this paper we merged the three approaches: We combined MALM with the stabilized-regularized path-following approach. We then put together the merit functions in \cite{MALM} and \cite{ForsgrenGill} to globalize the combined approach.

Eventually, we discussed on the issue of loss of centrality when the Augmented Lagrangien multipliers are updated; a problem that is not addressed in \cite{Renke,ALM_IPM1}. Utilizing the polynomially efficient solution of \eqref{eqn:KKT_LP}, we were able to set up an iteration scheme that is at least polynomially efficient in the special case where $f,c$ are linear. We believe that this is a promising feature of our method, that can make it particularly attractive for the case of solving \eqref{eqn:PBNLP} when $f,c$ are nonlinear.

\FloatBarrier

\bibliography{ALM_IPM_bib}
\bibliographystyle{plain}

\end{document}